\documentclass[12pt,reqno,a4wide]{amsart}

\usepackage{tikz} 
\usepackage{calrsfs}
\usepackage{mathrsfs}

\usepackage{array}
\usepackage{hyperref}

\usetikzlibrary{decorations.pathreplacing}
\usetikzlibrary{fit}

\usepackage{pgfplots}

\allowdisplaybreaks

\newtheorem{theorem}{Theorem}

\catcode`,\active

\catcode`\,12

\begin{document}
	\bibliographystyle{plain}
	
	%
	%
	
	\title[Motzkin paths with two variants of level steps on odd levels]
	{Motzkin paths with two variants of level steps on odd levels --- a kernel method approach}

	\author[H. Prodinger ]{Helmut Prodinger }
	\address{Department of Mathematics, University of Stellenbosch 7602, Stellenbosch, South Africa
		and
		NITheCS (National Institute for
		Theoretical and Computational Sciences), South Africa.}
	\email{warrenham33@gmail.com}


	\begin{abstract}
The sequence  A176677 in the Encyclopedia of Integer Sequences enumerates Motzkin paths where two types of horizontal steps may occur, but only on odd indexed levels.
We show how to perform the enumeration, also dealing with partial such Motzkin paths leading to a particular level or to any level (open paths). The method is the kernel method where functional equations are manipulated in a suitable way.  The coefficients of sequence  A176677 satisfy a holonomic recursion that was recently discussed on the arxiv. We show how this can be  established in an (almost) automatic fashion. Eventually we switch the roles of `odd' and `even'. One could also allow more versions of horizontal steps but we leave this to the interested readers.
		
	\end{abstract}

	\maketitle

	\section{Motzkin paths with variations on the possible level steps}

Motzkin paths have up steps $(1,1)$, down steps $(1,-1)$ and level steps $(1,0)$. Like Dyck paths, they do not go below the $x$-axis and return eventually to the $x$-axis.
Partial paths that do not necessarily return to the $x$-axis are also of interest. If there are two versions of level steps allowed, but only on an odd level, the enumeration leads to 
the sequence  A176677 in \cite{OEIS}\footnote{The encyclopedia provides 4 different references for this particular sequence}. We show how to compute this using the kernel method, as explained in \cite{Prodingerkernel} and many other papers.

We use an automaton to describe the Motzkin paths in question. We use two levels (`even' on top, `odd' at the bottom); as usual, there are loops, but for the bottom layer there are loops of two possible level steps.
\begin{figure}[h]\label{pso1}
	\begin{tikzpicture}[line width=1pt,scale=1.4]
		
		\draw (0,0) circle (0.05cm);
		\fill (0,0) circle (0.1cm);
		\draw (8.5,-0.5)node  {$\dots$};
		\draw (0,-1) circle (0.03cm);
		
		\foreach \t in{0,1,2,3,4,5,6,7}
		{
			\draw[red,latex-](\t+1,-1) to  [in=-10,out=100](\t,0);
			\draw[red,latex-](\t,-1) to  [in=210,out=60](\t+1,0);
			\draw[blue,-latex](\t,-1) to  [in=250,out=30](\t+1,0);
				\draw[blue,-latex](\t+1,-1) to  [in=300,out=150](\t,0);
			
			\draw (\t+1,0)circle (0.03cm);
			\draw (\t+1,-1)circle (0.03cm);
			\draw[-latex](\t,0) to [in=-10,out=10](\t-0.06,0.5);
			\draw[](\t,0) to [out=170,in=-170](\t,0.5);
			
			\draw[brown, ultra thick,-latex](\t,-1) to [in=10,out=-10](\t-0.06,-1.5);
			\draw[brown, ultra thick](\t,-1) to [out=-170,in=170](\t-0.00,-1.5);
		}
	\foreach \t in{8}
	{
		\draw (\t,0)circle (0.03cm);
		\draw[-latex](\t,0) to [in=-10,out=10](\t-0.06,0.5);
		\draw[](\t,0) to [out=170,in=-170](\t,0.5);
		\draw[brown,-latex, ultra thick](\t,-1) to [in=10,out=-10](\t-0.06,-1.5);
		\draw[brown, ultra thick](\t,-1) to [out=-170,in=170](\t-0.00,-1.5);
	}

	\end{tikzpicture}
	\caption{ Distinguishing even and odd levels (even levels on top). On odd levels, two variants of level steps are possible (depicted in brown). Odd numbered states on top and even numbered states at the bottom cannot be reached.}
\end{figure}
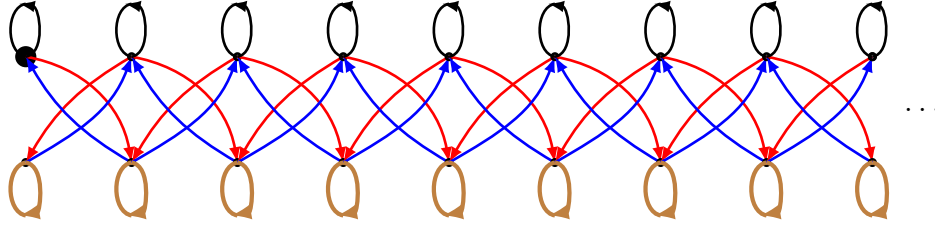

Since not all the states are accessible, one might simplify the automaton (Figure~2):
\begin{figure}[h]
	\begin{tikzpicture}[line width=1pt,scale=1.4]
		
		\draw (0,0) circle (0.05cm);
		\fill (0,0) circle (0.1cm);
		\draw (8.5,-.5)node  {$\dots$};
		\draw (8,0)circle (0.03cm);
		\foreach \t in{0,2,4,6}
		{
			\draw[red,latex-](\t+1,0-1) to  [out=100,in=0](\t,0);
			\draw[blue,latex-](\t,0) to  [in=180,out=-90](\t+1,0-1);
			
			\draw[red,latex-](\t+1,0-1) to  [out=80,in=180](\t+2,0);
				\draw[blue,-latex](\t+1,0-1) to  [out=0,in=-90](\t+2,0);
				
			\draw (\t,0)circle (0.03cm);
			\draw (\t+1,-1)circle (0.03cm);
			\draw[-latex](\t,0) to [in=-10,out=10](\t-0.06,0.5);
			\draw[](\t,0) to [out=170,in=-170](\t,0.5);
			
		}
	
	\foreach \t in{0,2,4,6,8}
	{
		\draw[-latex](\t,0) to [in=-10,out=10](\t-0.06,0.5);
		\draw[](\t,0) to [out=170,in=-170](\t,0.5);
		
	}
\foreach \t in{1,3,5,7}
{
	\draw[brown,-latex, ultra thick](\t,-1) to [in=10,out=-10](\t-0.06,-1.5);
	\draw[brown, ultra thick](\t,-1) to [out=-170,in=170](\t-0.00,-1.5);
}

	\end{tikzpicture}
\label{pso2}
	\caption{Simplified version}
\end{figure}
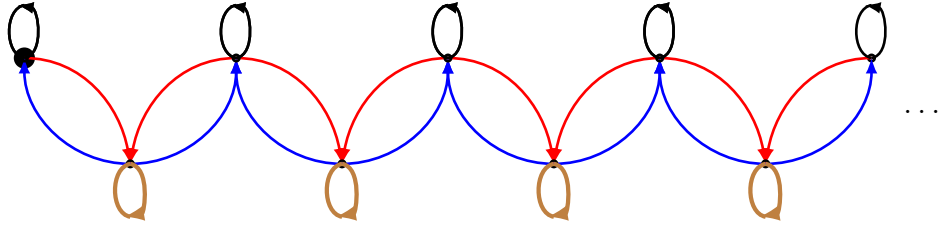

However, we find it more straight forward to stay with the two layers of states even though not all of them can be reached. We use generating functions $f_i(z)$ and $g_i(z)$
leading to level $i$ on top resp.\ bottom layer, and derive recursions according to what the last step was:

Notice that $f_1=f_3=f_5\dots=0$ and $g_0=g_2=g_4\dots=0$.
\begin{align*}
f_0&=1+zf_0+zg_1,\\
f_i&=zf_i+zg_{i+1}+zg_{i-1},\ i\ge1\\
g_i&=2zg_i+zf_{i+1}+zf_{i-1},\ i\ge1.
\end{align*}
Then we need bivariate generating functions as well:
\begin{equation*}
F(z,u)=\sum_{i\ge0}u^if_i(z),\quad G(z,u)=\sum_{i\ge0}u^ig_i(z);
\end{equation*}
$F(z,u)$ has only even powers of $u$, and $G(z,u)$ has only odd powers of $u$.
By summing up the recursions one finds
\begin{align*}
F(z,u)&=\sum_{i\ge0}u^if_i=1+zf_0+zg_1+\sum_{i\ge1}u^i[zf_i+zg_{i+1}+zg_{i-1}]\\
&=1 +\sum_{i\ge0}u^izf_i+\sum_{i\ge0}u^izg_{i+1}+\sum_{i\ge1}u^izg_{i-1}\\
&=1 +zF(z,u) +\frac zu[G(z,u)-g_0]+zu G(z,u).
\end{align*}
Notice again that $g_0=0$;
\begin{align*}
	G(z,u)&=\sum_{i\ge0}u^ig_i=\sum_{i\ge1}u^i[2zg_i+zf_{i+1}+zf_{i-1}]\\*
	&=2zG(z,u)+\frac zu[F(z,u)-f_0]+zu F(z,u).
\end{align*}
For simplicity, we will occasionally write $F(u)=F(z,u)$, $G(u)=G(z,u)$, $F(0)=f_0$, $G(0)=g_0=0$.
Solving the system,
\begin{equation*}
F(u)=\frac {2 z{u}^{2}-{u}^{2}+{z}^{2}{u}^{2}f_0
	+f_0 {z}^{2}}{{z}^{2}+{z}^{2}{u}^{4}-{u}^{2}+3 z{u}^{2}},
\end{equation*}
\begin{equation*}
G(u)=-{\frac {z ( u+uzf_0+{u}^{3}-uf_0 ) }{{z}^{2}+{z}^{2}{u}^{4}-{u}^{2}+3 z{u}^{2}}}.
\end{equation*}
We need to factor the denominator and abbreviate
\begin{equation*}
W:=\sqrt{(1-z)(1-2z)(1-3z-2z^2)};
\end{equation*}
\begin{align*}
s_1&=\frac{\sqrt{2}\sqrt{1-3z+W}}{2z},\quad s_2=-\frac{\sqrt{2}\sqrt{1-3z+W}}{2z},\\
 \ s_3&=\frac{\sqrt{2}\sqrt{1-3z-W}}{2z},\quad s_4=-\frac{\sqrt{2}\sqrt{1-3z-W}}{2z}.
\end{align*}
Then
\begin{equation*}
{z}^{2}+{z}^{2}{u}^{4}-{u}^{2}+3 z{u}^{2}=z^2(u-s_1)(u-s_2)(u-s_3)(u-s_4).
\end{equation*}
Plugging in $u=0$ leads to a void equation, but
the factors $u-s_3$ and $u-s_4$ are `bad' (in the spirit of the kernel method \cite{Prodingerkernel}, as they do not have a power series expansion around $(u,z)\sim(0,0)$) and must cancel out. Then (the cancellation can be done with Maple, regardless what $s_3$ and $s_4$ are; Maple has the \textsf{quo} command  for that) really are
\begin{align*}
F(u)&=\frac{2z-1+f_0z^2}{z^2(u-s_1)(u-s_2)},\\
G(u)&=\frac{-zu}{z^2(u-s_1)(u-s_2)}.
\end{align*}
Maple also has the \textsf{rem} command that we do not use here. Together, these two functions \textsf{quo} and \textsf{rem} provide the division with remainder. 
From this we find, by setting $u=0$
\begin{align*}
f_0&=\frac{1-3z+2z^2-W}{2(1-z)z^2}\\
&=1+z+2z^2+5z^3+14z^4+41z^5+123z^6+375z^7+1158z^8+3615z^9+\dots,
\end{align*}
which is the sequence A176677 in \cite{OEIS}.
Note that
\begin{equation*}
(u-s_1)(u-s_2)=u^2-s_1^2.
\end{equation*}
Then
\begin{equation*}
	G(u)=\frac{-u}{z(u-s_1)(u-s_2)}=\frac{-u}{z(u^2-s_1^2)}=\frac{u}{zs_1^2(1-u^2/s_1^2)}
\end{equation*}
and
\begin{equation*}
	[u^{2k+1}]G(u)=[u^{2k}]\frac{1}{zs_1^2(1-u^2/s_1^2)}=\frac1{zs_1^{2k+2}}.
\end{equation*}
Further
\begin{equation*}
F(u)=\frac{2z-1+f_0z^2}{z^2(u-s_1)(u-s_2)}=\frac{2z-1+f_0z^2}{z^2}\frac1{u^2-s_1^2}=\frac{2z-1+f_0z^2}{z^2}\frac{-1}{s_1^2}\frac1{1-u^2/s_1^2}
\end{equation*}
and therefore
\begin{equation*}
[u^{2k}]	F(u)=\frac{2z-1+f_0z^2}{z^2}\frac{-1}{s_1^2}[u^{2k}]\frac1{1-u^2/s_1^2}=\frac{2z-1+f_0z^2}{z^2}\frac{-1}{s_1^{2k+2}}
\end{equation*}
A direct computation yields
\begin{equation*}
	f_0=\frac1{z^2}-\frac2z-\frac{1+s_1^2}{1-z}
	\end{equation*}
and
\begin{equation*}
-\frac{2z-1+f_0z^2}{z^2}=	\frac{1+s_1^2}{1-z} .
\end{equation*}
Therefore
\begin{equation*}
[u^{2k}]F(u)=\frac{1+s_1^2}{1-z}\frac1{s_1^{2k+2}}.
\end{equation*}

Note that $F(1)+G(1)$ is the generating function of all \emph{open} Motzkin paths (variations as discussed in this note).
\begin{align*}
F(1)+G(1)=1+2z+6z^2+19z^3+62z^4+205z^5+684z^6+2298z^7+7764z^8+26355z^9+\dots
\end{align*}
This sequence of coefficients $1,2,6,19,62,205,684,2298,7764,26355\dots$ is not in the \cite{OEIS}.

We summarize our findings:
\begin{theorem}
Let 
\begin{equation*}
	s_1^2=\frac{{1-3z+\sqrt{(1-z)(1-2z)(1-3z-2z^2)}}}{2z^2},
\end{equation*}
then the generating function of Motzkin paths with 2 variations of level steps on odd levels, leading to level $2k$ is
\begin{equation*}
\frac{1+s_1^2}{1-z}\frac1{s_1^{2k+2}}
\end{equation*}
and to  level $2k+1$ is
\begin{equation*}
\frac1{zs_1^{2k+2}}.
\end{equation*}
The generating function of open Motzkin paths is
\begin{align*}
	F(1)+G(1)=1+2z+6z^2+19z^3+62z^4+205z^5+684z^6+2298z^7+7764z^8+26355z^9 +\dots;
\end{align*}
an analytic expression for $F(1)+G(1)$ is available from the text.
\end{theorem}

\section{Automatic derivation of  a holonomic recurrence for the Motzkin numbers in question}

We want to derive  a recursion for the coefficients of the generating function
\begin{align*}
	f_0&=\frac{1-3z+2z^2-\sqrt{(1-z)(1-2z)(1-3z-2z^2)}}{2(1-z)z^2};
\end{align*}
the software \textsf{gun} \cite{gfun} contains all the necessary tools for that. The following steps will be performed.
\begin{itemize}
	\item Find an algebraic equation that $f_0$ satisfies.
	\item Translate the algebraic equation into a differential equation.
	\item Translate the differential equation into a recurrence for the coefficients.
\end{itemize}

\begin{equation*}
>\textsf{algfuntoalgeq}(f_0, y(z));
\end{equation*}
the result is
\begin{equation*}
z^2(-1+z)y^2+(1-z)(1-2z)y+2z-1=0.
\end{equation*}
The algebraic equation can be translated into a differential equation
\begin{equation*}
>\textsf {algeqtodiffeq}(y^2z^3-y^2z^2+2yz^2-3yz+2z+y-1=0, y(z));
\end{equation*}
the result is
\begin{equation*}
6z^2-7z+2+(10z-14z^2+4z^3+4z^4-2)y(z)+(4z^5-9z^3-z+6z^2)y'(z)=0.
\end{equation*}
One can get a homogenous differential equation as well:
\begin{equation*}
> \textsf{diffeqtohomdiffeq}(6z^2-7z+2+(10z-14z^2+4z^3+4z^4-2)y(z)+(4z^5-9z^3-z+6z^2)y'(z)=0);
\end{equation*}
the result is
\begin{gather*}
2(24z^4-6z^3-18z^2+13z-3)(z-1)y(z)\\{}+2(2z-1)(24z^5-17z^4-38z^3+54z^2-23z+3)y'(z)\\ {}+z(z-1)(3z-2)(2z^2+3z-1)(2z-1)^2y{''}(z)=0.
\end{gather*}
Now we translate into a recursion for the (Taylor) coefficients
\begin{equation*}
>\textsf{diffeqtorec}(6z^2-7z+2+(10z-14z^2+4z^3+4z^4-2)y(z)+(4z^5-9z^3-z+6z^2)y'(z), y(z), m(n));
\end{equation*}
the result is the recursion
\begin{gather*}
4(n+1)m(n)+4m(n+1)-(9n+32)m(n+2)\\+2(3n+14)m(n+3)-(n+6)m(n+4)=0,\\ m(0)=1,\ m(1)=1,\ m(2)=2,\ m(3)=5.
\end{gather*}
This recursion of the coefficients 
\begin{equation*}
f_0(z)=\sum_{n\ge0}m(n)z^n
\end{equation*}
is equivalent to the one given in \cite{tongniu} but derived in a (semi-)automatic fashion. \textsf{gfun} has several additional features of practical interest \cite{gfun}.

\section{Switching the roles of even resp. odd levels }

As the section title says, we switch to roles of `even' resp.\ `odd'. The notation stays the same, and there are just a few modifications necessary.
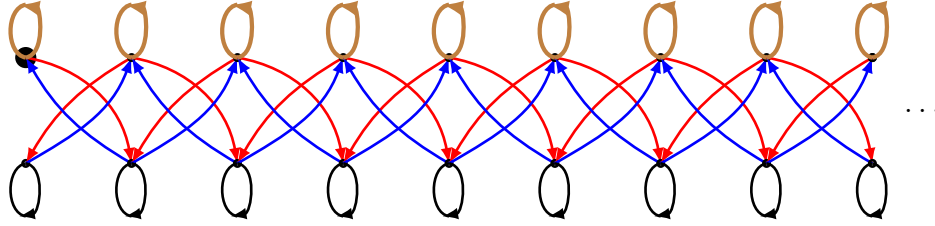
\begin{figure}[h]\label{pso3}
	\begin{tikzpicture}[line width=1pt,scale=1.4]
		
		\draw (0,0) circle (0.05cm);
		\fill (0,0) circle (0.1cm);
		\draw (8.5,-0.5)node  {$\dots$};
		\draw (0,-1) circle (0.03cm);
		
		\foreach \t in{0,1,2,3,4,5,6,7}
		{
			\draw[red,latex-](\t+1,-1) to  [in=-10,out=100](\t,0);
			\draw[red,latex-](\t,-1) to  [in=210,out=60](\t+1,0);
			\draw[blue,-latex](\t,-1) to  [in=250,out=30](\t+1,0);
			\draw[blue,-latex](\t+1,-1) to  [in=300,out=150](\t,0);
			
			\draw (\t+1,0)circle (0.03cm);
			\draw (\t+1,-1)circle (0.03cm);
			\draw[brown, ultra thick,-latex](\t,0) to [in=-10,out=10](\t-0.06,0.5);
			\draw[brown, ultra thick](\t,0) to [out=170,in=-170](\t,0.5);
			
			\draw[-latex](\t,-1) to [in=10,out=-10](\t-0.06,-1.5);
			\draw[](\t,-1) to [out=-170,in=170](\t-0.00,-1.5);
		}
		\foreach \t in{8}
		{
			\draw (\t,0)circle (0.03cm);
			\draw[brown, ultra thick,-latex](\t,0) to [in=-10,out=10](\t-0.06,0.5);
			\draw[brown, ultra thick](\t,0) to [out=170,in=-170](\t,0.5);
			\draw[-latex](\t,-1) to [in=10,out=-10](\t-0.06,-1.5);
			\draw[](\t,-1) to [out=-170,in=170](\t-0.00,-1.5);
		}

	\end{tikzpicture}
	\caption{Brown loops refer to two variants of level steps, black loops just to one variant.}
\end{figure}
Then
\begin{align*}
	F(z,u)	&=1 +2zF(z,u) +\frac zuG(z,u)+zu G(z,u),\\
	G(z,u)	&=zG(z,u)+\frac zu[F(z,u)-f_0]+zu F(z,u).
\end{align*}
Solving the system,
\begin{gather*}
	F(u)=\frac {z^2u^2f_0+z^2f_0+zu^2-u^2}{{z}^{2}+{z}^{2}{u}^{4}-{u}^{2}+3 z{u}^{2}},\\*
	G(u)=-\frac {zu(1-f_0+2zf_0+u^2)}{{z}^{2}+{z}^{2}{u}^{4}-{u}^{2}+3 z{u}^{2}}.
\end{gather*}
The denominators are the same as before. Consequently, dividing out the two bad factors $(u-s_3)(u-s_4)$,
\begin{align*}
	f_0&=\frac{1-3 z+2 z^2-W}{2(1-2z)z^2}\\
	&=1+2 z+5 z^2+13 z^3+35 z^4+97 z^5+276 z^6+804 z^7+2391 z^8+\dots,
\end{align*}
a sequence that is not in the OEIS \cite{OEIS}. Furthermore
\begin{align*}
F(u)&=\frac{z^2f_0+z-1}{z^2(u-s_1)(u-s_2)},\\
G(u)&=\frac{-u}{z(u-s_1)(u-s_2)}.
\end{align*}
The function $G(u)$ has not changed from before (which makes sense once one thinks about it). However,
\begin{equation*}
F(u)=\frac{\frac{1-3 z+2 z^2-W}{2(1-2z)}+z-1}{z^2(u^2-s_1^2)}=-\frac{1+s_1^2}{(1-2z)(u^2-s_1^2)}=\frac{1+s_1^2}{(1-2z)s_1^2(1-u^2/s_1^2)}
\end{equation*}
and
\begin{equation*}
[u^{2k}]F(u)=f_{2k}=\frac{1+s_1^2}{(1-2z)s_1^{2k+2}}.
\end{equation*}
Again, $f_{2k}$ and $g_{2k+1}$ can be expressed by the function $s_1^2$.

\end{document}